УДК 513.83; 517.5

Yu.B.Zelinskii, I.Yu.Vyhovs'ka, M.V.Stefanchuk, Generalized convex sets and shadow's problem .- Ukr.math.journ.

The problem of shadow is solved. It is equivalent to condition for point is in generalized convex hull of a family of compact sets.

Ю.Б.Зелінський, І.Ю.Виговська, М.В.Стефанчук, Узагальнено опуклі множини та задача про тінь.- Укр.матем.ж.

В роботі отримано повний розв'язок проблеми про тінь, що еквівалентно знаходженню умов належності точки узагальнено опуклій оболонці сім'ї компактних множин.

УДК 513.83; 517.5

Ю.Б.Зелинский, И.Ю.Выговская, М.В.Стефанчук (Ин-т математики НАН Украины, Киев)

## Обобщённо выпуклые множества и задача о тени

Главная цель работы – решение задачи о тени, которую можно рассматривать как нахождение условий обеспечивающих принадлежность точки обобщенно выпуклой оболочке некоторого семейства множеств.

**Определение 1.** Скажем, что множество $E \subset \mathbb{R}^n$ $m$-выпукло относительно точки $x \in \mathbb{R}^n \setminus E$, если найдется $m$-мерная плоскость $L$, такая что $x \in L$ и $L \cap E = \varnothing$.

**Определение 2.** Скажем, что множество $E \subset \mathbb{R}^n$ $m$-выпукло, если оно $m$-выпукло относительно каждой точки $x \in \mathbb{R}^n \setminus E$.

Легко убедиться, что оба приведенные определения удовлетворяют известной аксиоме выпуклости: пересечение каждого подсемейства таких множеств тоже удовлетворяет определению. Для произвольного множества $E \subset \mathbb{R}^n$ мы можем рассматривать минимальное $m$–выпуклое множество, содержащее $E$, и назвать его $m$-оболочкой множества $E$.

Как частный случай принадлежности точки 1-оболочке объединения некоторого набора шаров можно привести следующую задачу о тени, рассмотренную Г.Худайбергановым [1-3].

**Задача** (о тени). Какое минимальное число попарно непересекающихся замкнутых шаров с центрами на сфере $S^{n-1}$ и радиуса меньшего от радиуса сферы достаточно чтобы любая прямая, проходящая через центр сферы, пересекала хотя бы один из этих шаров?

Другими словами эту задачу можно переформулировать так. Сколько замкнутых шаров радиуса меньшего от радиуса сферы с центрами на сфере (минимальное количество) обеспечит принадлежность центра сферы 1-оболочке семейства шаров?

Если в сферу вписать правильный $n$–мерный симплекс и разместить шары радиуса равного половине длины ребра симплекса в его вершинах, то очевидно, что эта система шаров создает тень для центра сферы. Однако при этом мы нарушим одно условие – шары попарно касаются друг друга. Пусть $a$ половина длины ребра правильного симплекса. Рассмотрим семейство из $n+1$ шаров радиусов $a+\varepsilon$, $a-\varepsilon/2$, $a-\varepsilon/2^2$,…, $a-\varepsilon/2^n$, соответственно, для достаточно малого числа $\varepsilon$. Разместим эти шары так, чтобы они попарно касались друг друга, а их центры образовывали симплекс мало отличающийся от правильного. Через центры этих шаров проходит единственная сфера, центр которой принадлежит 1-оболочке семейства шаров. Внутренности этого семейства шаров образуют семейство из $n+1$ открытого шара, для которого центр сферы принадлежит 1-оболочке этого семейства. Если же исходные замкнутые шары немного уменьшить, то в силу непрерывности, очевидно, что $n+1$ замкнутого шара достаточно для создания тени.

**Лемма 1.** Если множество $E = \bigcup_{i=1}^{n-1} E_i \subset \mathbb{R}^n$ представляет собой объединение из $n$-1 выпуклого множества, то $E$ – 1-выпуклое множество.

Доказательство. В силу выпуклости каждого множества $E_i$ для произвольной точки $x \in \mathbb{R}^n \setminus E$ существует гиперплоскость $L_i$, содержащая эту точку, которая не пересекает множество $E$. Пересечение этих гиперплоскостей $L = \bigcap_{i=1}^{n-1} L_i$ содержит искомую прямую.

Из этой леммы следует, что произвольной совокупности из $n$-1 шаров для создания тени мало. Поэтому точное значение необходимого количества шаров $n$ или $n+1$.

Аналогично легко доказать следующее утверждение.

**Следствие 1.** Если множество $E = \bigcup_{i=1}^{m-1} E_i \subset \mathbb{R}^n$ представляет собой объединение из $m$-1 выпуклого множества, где $m < n$, то $E$ – ($n$-$m$)-выпуклое множество.

Задача о тени была решена Г.Худайбергановым для *n* = 2 (показано, что двух шаров достаточно). Здесь же было предложено решение при *n* > 2, которое оказалось ошибочным и насколько известно авторам точное значение количества шаров остается открытой проблемой. Дальнейшие рассуждения позволят дать полный ответ на эту проблему.

Приведем здесь другое решение задачи о тени для *n* = 2, использующее непрерывность изменения прямых, и дающее некоторые цифровые оценки. Исследуем 1-оболочку объединения двух шаров $K_1$, $K_2$ в случае *n* = 2. Для этой пары шаров есть касательная прямая $l_1$, которая разделяет шары, и касательная прямая $l_2$ для которой оба шара находятся в одной полуплоскости, задаваемой прямой. Мы рассмотрим предельный в некотором смысле случай, когда шары касаются друг друга и кроме этого прямая $l_2$ проходит через центр сферы $S^1$ (окружности). При этом прямая $l_1$ естественно проходит через точку касания шаров. Пусть для радиусов шаров справедливо неравенство $1 \geq r_1 \geq r_2$, соответственно. Используя теорему Пифагора легко показать, что точки касания шаров вырезают из прямой $l_2$ отрезок длиной $2\sqrt{r_1 r_2}$, а расстояние от центра окружности до прямой $l_1$ равно $(r_1 - r_2)/2$. Узнаем, когда последнее расстояние можно сделать максимальным. Для этого положим $r_1 = 1$ и найдем значение $r_2$, используя то, что теперь точка касания шара $K_1$ с прямой $l_2$ совпадает с центром окружности. Получим равенство

$$2\sqrt{r_2} = \sqrt{1 - r_2^2}.$$

Дальше все сводится к квадратному уравнению

$$r_2^2 + 4r_2 - 1 = 0,$$

один из корней которого $r_2 = \sqrt{5} - 2$ показывает, что максимально возможное отклонение прямой $l_1$ от центра окружности равно $(3 - \sqrt{5})/2$. Теперь становиться ясным, почему мы выбрали прямую $l_2$ проходящей через центр окружности. Увеличение радиуса шара $K_2$ уменьшило бы искомое расстояние. Однако выбранные шары пока не удовлетворяют условию задачи о тени. Радиус большего шара равен 1 и шары касаются друг друга. Но, пользуясь, как и выше непрерывностью изменения прямой мы можем чуть уменьшить радиус большего шара и раздвинуть немного центры шаров, так что при этом прямая $l_1$ все ещё не сможет пройти через центр окружности. Кстати при этом шары, можно выбрать открытыми. Очевидно, что мы можем, как угодно близко приблизится к найденной константе $(3-\sqrt{5})/2$. Отсюда справедлива теорема.

**Теорема 1.** Существует два замкнутых (открытых) шара с центрами на единичной окружности и радиуса меньше 1, которые обеспечивают принадлежность центра окружности 1-оболочке семейства шаров.

Совершенно аналогично доказывается следующее утверждение.

**Следствие 2.** Существует два замкнутых (открытых) шара с центрами на сфере $S^{n-1}$ и радиуса меньшего от радиуса сферы, которые обеспечивают принадлежность центра сферы ( $n$- 1)-оболочке семейства шаров.

Покажем, что в случае двумерной сферы трех шаров недостаточно. Точки пространства будем обозначать координатами ($x$, $y$, $z$). Не нарушая общности, предположим, что сфера с центром в начале координат имеет радиус 1 и, что трех открытых шаров достаточно для создания тени в центре сферы. Очевидно, что эти шары должны иметь разные радиусы, потому что при равенстве радиусов даже при касании шаров через точку касания двух шаров проходит двумерная плоскость, которая может пересечь только третий шар. Поэтому в этой плоскости через центр сферы, согласно геометрической форме теоремы Хана-Банаха [4], проходит прямая не пересекающая ни одного шара. Предположим, что имеют место неравенства для радиусов шаров 1≥$r_1$>$r_2$>$r_3$. Мы можем считать, что шары попарно касаются друг друга, иначе мы могли бы их увеличить с тем же эффектом для тени. Расположим центр шара максимального радиуса в точке (0,0,1). Проведем двумерную плоскость $L$ через центры шаров $B_1$, $B_2$ и начало координат (Рис.1). Каждый из этих шаров порождает круговой конус с центром в начале координат такой, что произвольная прямая, лежащая внутри конуса, пересекает этот шар. Этот конус пересекает сферу по двум окружностям, которые можно задать парой параллельных плоскостей. Им в плоскости $L$ соответствуют две прямые $GB$ и $CA$. Полоса между этими плоскостями вырезает из сферы часть тех точек каждая прямая через центр сферы и такую точку не пересекает соответственный шар. Аналогично для второго шара получим две прямые $CG$ и $AB$ в плоскости $L$.

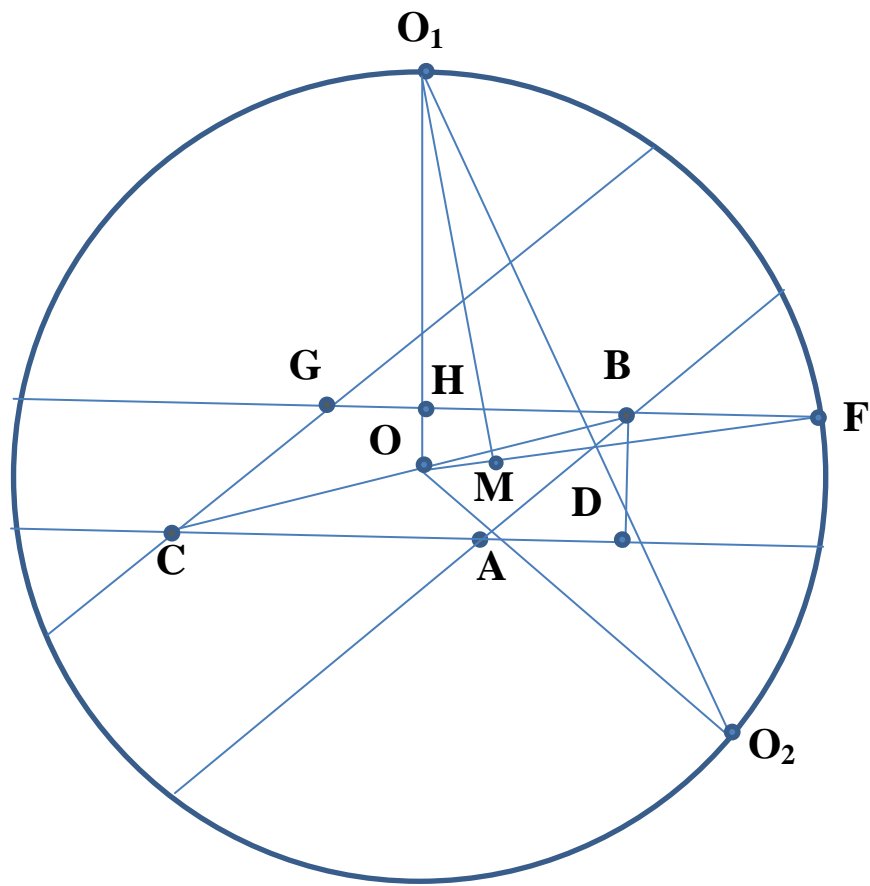

Рис.1

Пересечение двух полос, соответствующих двум шарам, представляет собой цилиндр, в основании которого лежит параллелограмм $ABGC$. Этот цилиндр вырезает на сфере множество тех точек каждая прямая через центр сферы и такую точку не пересекает оба шара. Теперь очевидно, что радиус третьего шара необходимого для создания тени не может быть меньше половины диагонали $BC$ этого параллелограмма. Прямая $OF$ касается шара $B_1$ в точке *M*, поэтому из равенства треугольников $OO_1M$ и $OHF$ следует, что длина *HF* равна $r_1$. Отсюда получим, что расстояние между прямыми *GB* и *CA* равно $|BD| = 2\sqrt{1-r_1^2}$. Поскольку шары $B_1$ и $B_2$ касаются, то $|O_1O_2| = r_1 + r_2$. Следовательно $\angle O_1OO_2 = 2\arcsin((r_1 + r_2)/2)$. Отрезок $OO_1$ перпендикулярен к прямой *GB*, а отрезок $OO_2$ перпендикулярен к прямой *AB*. Отсюда найдем стороны параллелограмма *ABGC*. Имеем

$$|AB| = 2\sqrt{1-r_1^2} / \sin(2\arcsin((r_1 + r_2)/2)).$$

Аналогично

$$|AC| = 2\sqrt{1-r_2^2} / \sin(2\arcsin((r_1 + r_2)/2)).$$

Теперь из теоремы косинусов получим, что

$$|OB| = \frac{\sqrt{2 - r_1^2 - r_2^2 - 2\sqrt{1-r_1^2}\sqrt{1-r_2^2}\cos(2\arcsin((r_1+r_2)/2)))}}{\sin(2\arcsin((r_1+r_2)/2)))}. \qquad (1)$$

Произведем числовые оценки. Поскольку из вписанных в окружность треугольников максимальный периметр имеет правильный треугольник, то сумма радиусов шаров не превосходит полупериметра правильного треугольника, вписанного в единичную окружность

$$r_1 + r_2 + r_3 \le 1.5\sqrt{3} \approx 2.598. \qquad (2)$$

Радиус шара $B_2$ не может быть меньше $\sqrt{2}/2$, иначе из-за неравенства $r_2 > r_3$ шары $B_2$ и $B_3$ не смогут обеспечить пересечение с ними произвольной прямой, проходящей через начало координат и лежащей в плоскости *xOy*, а шар $B_1$ с этой плоскостью не пересекается.

Используя программу Derive из (1) следует, что при радиусе $r_2 < 0.77$ радиус $r_3 > 0.77$. Следовательно, не выполнено неравенство $r_2 > r_3$. Если же $r_2 > 0.85$, то также использованием программы Derive получаем неравенство $r_1 + r_2 + r_3 > 2.6$, что противоречит (2). Для суммы радиусов шаров $r_1 + r_2$ получаем неравенства $1.54 < 2r_2 < r_1 + r_2 \le 1 + r_2 < 1.85$.

Дальше получим следующие оценки (Рис.2)

$$1{,}1858 \leq |O_1K| = (r_1 + r_2)^2/2 \leq 1{,}7113,$$

$$0{,}9826 \geq |KO_2| = |OL| = (r_1 + r_2)\sqrt{1 - (r_1 + r_2)^2/4} \geq 0{,}7029,$$

$$0{,}1858 \leq |OK| = |O_1K| - 1 = |LO_2| \leq 0{,}7113,$$

$$0{,}4654 \leq |NL| = \sqrt{|NO_2|^2 - |LO_2|^2} = \sqrt{|r_2|^2 - |LO_2|^2} \leq 0{,}5216.$$

Отрезок $|NL|$ равен радиусу окружности, по которой шар $B_2$ пересекает плоскость $xOy$.

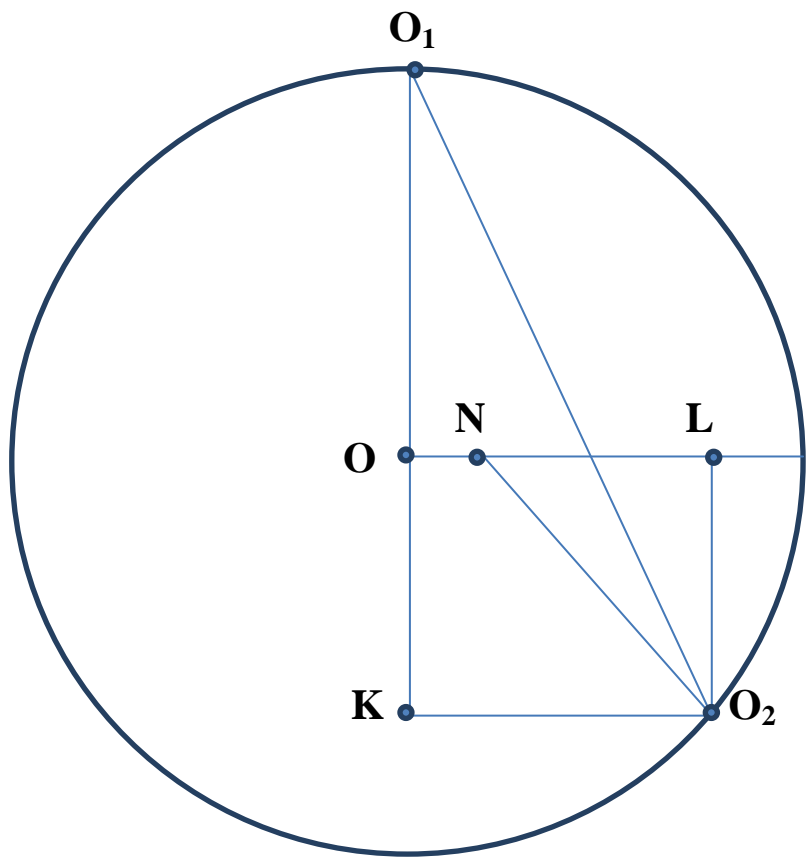

Рис.2

Точка $L$ - центр этой окружности. Отношение $|NL|/|OL|$ задает синус половины угла $\alpha$, под которым видна эта окружность из начала координат. Для $\sin(\alpha/2)$ имеем оценку сверху $\sin(\alpha/2) = |NL|/|OL| \leqslant 0{,}6621$, тогда arc $\sin(\alpha/2) \leqslant$ arc $\sin 0{,}6621 \leqslant 0{,}7236$. Поэтому угол $\alpha$ не превосходит 1,4472.

Аналогичные рассуждения в силу неравенства для радиусов показывают, что окружность, по которой с плоскостью $xOy$ пересекается шар $B_3$, так же видна из начала координат под углом не превышающим 1,4472. Эти две окружности суммарно закрывают угол, который не превышает 2,8945, что меньше развернутого угла размерности π. Поэтому в плоскости *xOy* существует прямая через начало координат не пересекающая ни одного из трех шаров. Следовательно, для создания тени в центре сферы при *n*=3 необходимо четыре шара. При *n* > 3 оценка получается применением математической индукции. Рассматриваем гиперплоскость через центр сферы, которая не пересекает один из шаров. Для создания тени в начале координат этой гиперплоскости, согласно предположению индукции, необходимо *n* шаров. Поэтому, прибавляя шар, который не пересекает выбранную гиперплоскость, получаем необходимость *n*+1-го шара. Получили следующее утверждение полностью решающее проблему тени.

**Теорема 2**. Для того чтобы центр (*n*-1)-сферы в *n*-мерном евклидовом пространстве при *n* > 2 принадлежал 1-оболочке семейства открытых (замкнутых) шаров радиуса не превышающего (меньшего) радиуса сферы и с центрами на сфере необходимо и достаточно *n*+1-го шара.

Рассмотрим более общие по отношению к предыдущим определениям объекты.

**Определение 3.** Скажем, что множество $E \subset \mathbb{R}^n$ *m*-полувыпукло относительно точки $x \in \mathbb{R}^n \setminus E$, если найдется *m*-мерная полуплоскость *P*, такая что $x \in P$ и $P \cap E = \varnothing$.

**Определение 4.** Скажем, что множество $E \subset \mathbb{R}^n$ *m*-полувыпукло, если оно *m*-полувыпукло относительно каждой точки $x \in \mathbb{R}^n \setminus E$.

Легко убедиться, что и эти определения удовлетворяют аксиоме выпуклости, и мы тоже можем строить *m*-полувыпуклые оболочки множеств согласно этим определениям.

Рассмотрим аналог задачи о тени для полувыпуклости. Какое минимальное число попарно непересекающихся замкнутых (открытых) шаров с центрами на сфере $S^{n-1}$ и радиуса меньшего от (не превышающего) радиуса сферы достаточно чтобы любой луч из центра сферы пересекал хотя бы один из этих шаров?

Задача проста в плоском случае $n = 2$. Если мы впишем в окружность остроугольный треугольник с неравными сторонами $a > b > c$, а в его вершинах разместим три замкнутых круга радиусов $p-a, p-b, p-c$, соответственно, где $p=(a+b+c)/2$, то очевидно, что полувыпуклая оболочка объединения этих кругов состоит из кругов и внутренности треугольника. Если центр окружности не принадлежит объединению кругов, то такая конструкция обеспечит тень и в этой точке. Теперь, как и выше, пользуясь непрерывностью, если чуть уменьшить радиусы кругов, то получим, что при $n = 2$ три замкнутых (открытых) круга решают задачу. Исследуем соотношение сторон треугольника, которые обеспечивают решение. Из неравенств $p-a < p-b < p-c$ следует, что радиус описанной окружности должен превышать $p-c$. Не нарушая общности, будем считать, что сторона $c=1$, и справедливы неравенства $a > b > 1$. Другие треугольники с нужным свойством получаются преобразованием подобия. Из формулы для радиуса описанной окружности, заменяя стороны треугольника переменными $x=a$, $y=b$, получим, что координаты нужных сторон должны находиться внутри криволинейного треугольника, две стороны которого прямые $x=y$, $y=1$, а третья кривая, заданная неявным уравнением,

$$x + y - 1 = \frac{2xy}{\sqrt{(x+y+1)(-x+y+1)(x-y+1)(x+y-1)}}.$$

Построением графика в Derive убеждаемся, что множество таких точек непустое. Следовательно, справедливо утверждение.

**Теорема 3.** Для того чтобы центр окружности $S^1 \subset \mathbb{R}^2$ принадлежал 1-полувыпуклой оболочке семейства открытых (замкнутых) кругов радиуса не превышающего (меньшего) радиуса окружности и с центрами на этой окружности необходимо и достаточно трех кругов.

С увеличением размерности задача усложняется. Покажем сначала, что существуют семейства выпуклых множеств, 1-полувыпуклая оболочка которых совпадает с таким семейством.

**Лемма 2.** Если множество $K = \bigcup_{i=1}^{n} K_i$, где все множества $K_i$ – выпуклые компакты, то $H^k(K) = 0 \quad k \geq n-1$ при $n > 1$, ($H^k(K)$ - группы когомологий компакта $K$ [5]).

• Доказательство. Докажем результат применением индукции. Как показано в [6] объединение двух выпуклых компактов не может быть носителем никакого коцикла в положительных размерностях. Следовательно теорема верна при $n = 2$. Предположим, что она верна при $m = n-1$ и докажем ее при $m = n$. Применим точную когомологическую последовательность Майера-Вьеториса [5] для триады

$$(\bigcup_{i=1}^{n} K_i, \bigcup_{i=1}^{n-1} K_i, K_n).$$

Выпишем три ее последовательных члена

$$H^j(\bigcup_{i=1}^{n-1}(K_i \cap K_n)) \to H^{j+1}(\bigcup_{i=1}^{n} K_i) \to H^{j+1}(\bigcup_{i=1}^{n-1} K_i) \oplus H^{j+1}(K_n).$$

Множества $(K_i \cap K_n)$ выпуклые. Поэтому в силу предположения индукции первый член последовательности нулевой при $j \geq n-2$. То же следует и для обоих слагаемых третьего члена последовательности. Теперь в силу точности последовательности получаем утверждение леммы.

**Теорема 4.** Каждое множество $K = \bigcup_{i=1}^{n} K_i$ в $\mathbb{R}^n$, где все множества $K_i$ – выпуклые компакты, является 1-полувыпуклым.

• Доказательство. Рассмотрим произвольную точку $x$ не лежащую в $K$. Выберем сферу $S^{n-1}$ с центром в точке $x$ так чтобы внутри шара, ограниченного этой сферой, точек $K$ не было. Для каждого компакта $K_i$ построим однополостный конус $conK_i$ с вершиной в точке $x$. Пусть множества $E_i$ заданы пересечениями $E_i =(con K_i) \cap S^{n-1}$. Рассмотрим их выпуклые оболочки [7] $F_i = conv E_i$. Каждое из построенных сейчас множеств выпукло и является подмножеством соответственного конуса. Объединение множеств $F_i$ не может содержать всю сферу $S^{n-1}$ иначе оно было бы носителем ненулевого коцикла, что невозможно согласно предыдущей лемме. Следовательно, на сфере $S^{n-1}$ найдется точка $y$ не принадлежащая объединению $\cup F_i$ и тогда луч выходящий из точки $x$ через точку $y$ не пересекает ни одного из множеств $K_i$. В силу произвольности выбора точки $x$ теорема доказана.

**Замечание.** Из множества ($n$-1)–мерных граней $n$–мерного симплекса легко составить множество, которое 1-полувыпуклым не будет.

Усложним задачу, наложив на множество дополнительные условия. Исследуем, когда семейство шаров с центрами на фиксированной сфере обеспечит принадлежность центра сферы 1-полувыпуклой оболочке семейства.

Пусть $S^2 \subset \mathbb{R}^3$ единичная сфера. Точки пространства будем обозначать координатами ($x$, $y$, $z$). Выберем два открытых шара единичного радиуса в точках (0,0,1) и (0,0,-1). Теперь лучи, которые не пересекают эти два шара, должны лежать в плоскости $xOy$. Открытый шар радиуса $\sqrt{2}-1$ в точке (1,0,0) касается заданных двух шаров и виден из начала координат в плоскости $xOy$ под углом $\alpha$, синус половины которого равен $\sqrt{2}-1$. Следовательно, $\alpha/2=\arc\sin(\sqrt{2}-1)=0.4271$, $\alpha=0.8542$. Поскольку этот угол помещается 7.35 раз в развернутом угле $2\pi$, то заполним окружность в плоскости $xOy$ четырьмя шарами радиуса $\sqrt{2}-1$ с центрами в точках (1,0,0), (0,1,0), (-1,0,0), (0,-1,0) соответственно. После этого между ними четыре шара радиусов $3-2\sqrt{2}$ с центрами в точках пересечения единичной окружности плоскости $xOy$ с биссектрисами координатных углов. Эти шары касаются двух соседних из предыдущих четырех. В силу разности радиусов соседних шаров с центрами в плоскости *xOy* этот набор из 10 шаров обеспечит принадлежность центра сферы 1-полувыпуклой оболочке их объединения. Как и выше чуть уменьшая радиусы шаров, видим, что существует набор замкнутых десяти шаров с теми же свойствами. Получаем следующее утверждение.

**Теорема 5**. Для того чтобы центр двумерной сферы в трехмерном евклидовом пространстве принадлежал 1-полувыпуклой оболочке семейства открытых (замкнутых) шаров радиуса не превышающего (меньшего) радиуса сферы и с центрами на сфере достаточно десяти шаров.

Вложением рассмотренного выше трехмерного пространства как линейного подпространства в $\mathbb{R}^n$ вместе с десятком шаров радиусов выбранных при доказательстве теоремы 5 получается следующая оценка для ($n$-2)-полувыпуклости.

**Следствие 3.** Для того чтобы центр ($n$-1)-мерной сферы в евклидовом пространстве $\mathbb{R}^n$ принадлежал ($n$-2)-полувыпуклой оболочке семейства открытых (замкнутых) шаров радиуса не превышающего (меньшего) радиуса сферы и с центрами на сфере достаточно десяти шаров.

К сожалению, предыдущие рассуждения не переносятся на более высокие размерности и не дают необходимых условий даже в трехмерном пространстве. Поэтому следующие вопросы остаются открытыми.

**Вопрос 1.** Какое минимальное количество шаров в трехмерном евклидовом пространстве обеспечит принадлежность центра сферы их 1-полувыпуклой оболочке?

В размерностях выше трех неясно даже существование конечного необходимого множества шаров.

**Вопрос 2.** Существует ли конечное количество шаров с перечисленными выше условиями в евклидовом пространстве $\mathbb{R}^n$, $n > 3$,

которое обеспечит принадлежность центра сферы их 1-полувыпуклой оболочке?

Если разрешить центрам шаров находиться на двух концентричных сферах, то, используя конструкцию для 1–выпуклости в $\mathbb{R}^n$, шары и радиус второй сферы получим гомотетией относительно центра первой сферы. Коэффициент гомотетии выберем с отрицательным знаком, так чтобы образ гомотетии не пересекался с исходным множеством. Очевидно, что центр сферы будет принадлежать 1-полувыпуклой оболочке шаров. Поэтому $2n+2$ шаров для этого достаточно.

# Литература.

Інститут математики НАНУ
Вул..Терещенківська,3
Київ
  e-mail: zel@imath.kiev.ua
vkirinata@gmail.com
stefanmv43@gmail.com